\documentclass[10.5pt]{article}
\usepackage[leqno]{amsmath}
\usepackage{amsfonts}
\usepackage{graphicx}

\newtheorem{condition**}{A*}
\newtheorem{condition***}{C*}
\newtheorem{condition*}{C}

\newtheorem{example}{Example}[section]
\newtheorem{proposition}{Proposition}[section]
\newtheorem{corollary}{Corollary}[section]
\newtheorem{definition}{Definition}[section]
\newtheorem{theorem}{Theorem}[section]
\newtheorem{lemma}{Lemma}[section]

\newtheorem{remark}{Remark}[section]

\numberwithin{equation}{section}

\newcommand{\dsp}{\displaystyle}

\DeclareMathOperator*{\esssup}{ess\, sup}
\def\ms{\medskip}
\def\q{\quad}
\def\qq{\qquad}
\def\3n{\negthinspace \negthinspace \negthinspace }
\def\2n{\negthinspace \negthinspace }
\def\1n{\negthinspace }
\def\ba{\begin{array}}
\def\ea{\end{array}}
\def\sqr#1#2{{\vcenter{\vbox{\hrule height.#2pt
              \hbox{\vrule width.#2pt height#1pt \kern#1pt \vrule width.#2pt}
              \hrule height.#2pt}}}}
\def\signed #1{{\unskip\nobreak\hfil\penalty50
              \hskip2em\hbox{}\nobreak\hfil#1
              \parfillskip=0pt \finalhyphendemerits=0 \par}}
\def\endpf{\signed {$\sqr69$}}

\def\dbE{\mathbb{E}}
\def\dbF{\mathbb{F}}

\def\dbP{\mathbb{P}}

\def\dbR{\mathbb{R}}
\def\dbS{\mathbb{S}}

\def\={\buildrel \triangle \over =}

\def\ds{\displaystyle}

\def\ns{\noalign{\ss}}

\def\d{\delta}
\def\e{\varepsilon}

\def\l{\lambda}

\def\si{\sigma}
\def\t{\tau}

\def\o{\omega}

\def\G{\Gamma}

\def\Th{\Theta}
\def\L{\Lambda}

\def\F{\Phi}
\def\Om{\Omega}

\def\cA{{\cal A}}

\def\cF{{\cal F}}

\def\cT{{\cal T}}
\def\cU{{\cal U}}

\def\ss{\smallskip}
\def\ms{\medskip}
\def\bs{\bigskip}
\def\q{\quad}
\def\qq{\qquad}
\def\hb{\hbox}

\def\liminf{\mathop{\underline{\rm lim}}}
\def\esssup{\mathop{\rm esssup}}

\def\da{\mathop{\downarrow}}
\def\Ra{\mathop{\Rightarrow}}

\def\lan{\mathop{\langle}}
\def\ran{\mathop{\rangle}}

\def\esssup{\mathop{\rm esssup}}

\def\wt{\widetilde}
\def\cd{\cdot}

\def\as{\hbox{\rm a.s.{ }}}

\def\({\Big (}
\def\){\Big )}
\def\[{\Big[}
\def\]{\Big]}
\def\bde{\begin{definition}}
\def\ede{\end{definition}}
\def\be{\begin{equation}}
\def\bel{\begin{equation}\label}
\def\ee{\end{equation}}
\def\bex{\begin{example}}
\def\eex{\end{example}}
\def\bt{\begin{theorem}}
\def\et{\end{theorem}}
\def\bc{\begin{corollary}}
\def\ec{\end{corollary}}
\def\bl{\begin{lemma}}
\def\el{\end{lemma}}
\def\bp{\begin{proposition}}
\def\ep{\end{proposition}}
\def\bas{\begin{assumption}}
\def\eas{\end{assumption}}
\def\br{\begin{remark}}
\def\er{\end{remark}}
\def\ba{\begin{array}}
\def\ea{\end{array}}
\def\ed{\end{document}}

\begin{document}

\title{A Mixed Linear Quadratic Optimal Control Problem \\
with a Controlled Time Horizon\footnote{This work is supported in
part by RGC Grants GRF521610 and GRF501010, and NSF Grant
DMS-1007514.}}

\author{Jianhui Huang$^a$, \quad Xun Li$^a$ \quad and \quad Jiongmin Yong$^b$ \\
{\footnotesize\textsl{$^a$Department of Applied Mathematics, The Hong Kong Polytechnic University, Hong Kong, China}} \\
{\footnotesize\textsl{$^b$Department of Mathematics, University of Central Florida, Orlando, FL 32816, USA}}}

\maketitle

\begin{abstract}

A \emph{mixed linear quadratic} (MLQ, for short) optimal control
problem is considered. The controlled stochastic system consists of
two diffusion processes which are in different time horizons. There
are two control actions: a standard control action $u(\cd)$ enters
the drift and diffusion coefficients of both state equations, and a
stopping time $\t$, a possible later time after the first part of
the state starts, at which the second part of the state is
initialized with initial condition depending on the first state. A
motivation of MLQ problem from a two-stage project management is
presented. It turns out that solving MLQ problem is equivalent to
sequentially solve a \emph{random-duration linear quadratic} (RLQ,
for short) problem and an \emph{optimal time} (OT, for short)
problem associated with Riccati equations. In particular, the
optimal cost functional can be represented via two coupled
stochastic Riccati equations. Some optimality conditions for MLQ
problem is therefore obtained using the equivalence among MLQ, RLQ
and OT problems.

\end{abstract}

\ms


\section{Preliminary and Problem Formulation}

Let $T>0$ be given and $(\Om,\cF,\dbP,\dbF)$ be a complete filtered
probability space on which a one dimensional standard Brownian
motion $W(\cd)$ is defined with $\dbF=\{\cF_t\}_{t\ge0}$ being its
natural filtration augmented by all the $\dbP$-null sets. We
consider the following stochastic controlled system:
\bel{state11}\left\{\ba{ll}
\ns\ds
dX_1(t)=\big[A_1(t)X_1(t)+B_1(t)u(t)\big]dt+\big[C_1(t)X_1(t)+D_1(t)u(t)
\big]dW(t),\q t\in[0,\t),\\
\ns\ds dX(t)=\big[A(t)X(t)+B(t)u(t)\big]dt+\big[C(t)X(t)+D(t)u(t)
\big]dW(t),\q t\in[\t,T],\\
\ns\ds X_1(0)=x_1,\qq X(\tau)=K(\t)X_1(\t-0), \ea\right. \ee
where
$A_1(\cd),B_1(\cd),C_1(\cd),D_1(\cd),A(\cd),B(\cd),C(\cd),D(\cd)$
and $K(\cd)$ are given matrix-valued functions of compatible sizes.
In the above, $X(\cd)=(X_1(\cd),X_2(\cd))$ is the state process,
taking values in $\dbR^n$, which is decomposed into two parts,
$X_i(\cd)$ is valued in $\dbR^{n_i}$ ($i=1,2$, $n_1+n_2=n$),
$u(\cd)$ is a (usual) control process taking values in some set
$U\subseteq\dbR^m$, and $\t$ is an $\dbF$-stopping time. From the
above, we see that the part $X_1(\cd)$ of the state process $X(\cd)$
starts to run from $x_1\in\dbR^{n_1}$ at $t=0$. The total system
will start to run at a later time $t=\t$, with the initial state
$X(\t)$ depending on $X_1(\t-0)$. Besides the usual control
$u(\cd)$, the stopping time $\t$ will also be taken as a control.
The above state equation can be interpreted as follows: we let
$X_1(\cd)$ represent the dynamics of some basic project whereas
$X_2(\cd)$, initialized at the time $\t$, represents an additional
or an auxiliary project. It is notable that the initial value of
$X_2(\cd)$ depends on $X_1(\t)$, the value of the first component of
the state at time $\t$. Some real examples are as follows.

\ms

\bf Example 1.1. \rm (Urban Planning) Let $X_1(t)$ denote some quantity of the
dynamic value of some basic infrastructure investment in urban
planning (for example, the transportation network, systematic
pollution protection, and so on) at time $t$, while $X_2(t)$ denotes
the quantity of the real-estate property in urban planning, at time
$t\ge\t$, where $\t\in(0,T)$ is the time moment at which some basic
infrastructure has been set. Note that the construction of basic
infrastructure will still be continued (although it will be less
intensive) after $\t$. It follows naturally the real estate property
should depend closely on $X_1(\cd)$.

\ms

\bf Example 1.2. \rm (Applied Technology) Let $X_1(\cd)$ represent
the capital investment of some high-tech company in the phrase of
primary research and development (R\&D) while let $X_2(\cd)$ denote
the capital investment in the phrase of technology marketing and
product promotion, etc. Of course, $X_2(\cd)$ will depend on the
competitive ability of the product which in turn depends on the
technology ability in basic research $X_1(\cd)$.

\ms

Note that if $A(\cd),B(\cd),C(\cd),D(\cd)$ are of the following
form:
$$A(\cd)=\left(\ba{cc}0&0\\0&A_2(\cd)\ea\right),\q
B(\cd)=\left(\ba{c}0\\B_2(\cd)\ea\right), \q
C(\cd)=\left(\ba{cc}0&0\\ 0&C_2(\cd)\ea\right),\q
D(\cd)=\left(\ba{c}0\\ D_2(\cd)\ea\right),$$
then on the time period $[\t,T]$, the part $X_1(\cd)$ will be
completely stopped and only the part $X_2(\cd)$ will be running.

\ms

Now we introduce the following quadratic cost functional:
\bel{cost}\ba{ll}
\ns\ds J(x_1;u(\cd),\t)={1\over2}\dbE\Big\{\int_0^\t\big[\lan Q_1(t)X_1(t),X_1(t)\ran+\lan R_1(t)u(t),u(t)\ran\big]dt+\lan G_1(\t)X_1(\t),X_1(\t)\ran\\
\ns\ds\qq\qq\qq\qq\q+\int_\t^T\big[\lan Q(t)X(t),X(t)\ran +\lan R(t)u(t),u(t)\ran\big]dt+\lan GX(T),X(T)\ran\Big\},\ea\ee
where $Q_1(\cd), R_1(\cd),Q(\cd),R(\cd)$, and $G_1(\cd)$ are
symmetric matrix-valued functions, and $G$ is a
symmetric matrix, of suitable sizes. Roughly speaking, our optimal
control problem is to minimize $J(x_1;u(\cd),\t)$ over the set of
all admissible controls $(u(\cd),\t)$. We now make our problem
formulation more precise.

\ms

For Euclidean space $\dbR^n$, we denote by $\lan\cd\,,\cd\ran$ its
inner product and $|\cd|$ the induced norm. Next, let $\dbR^{m\times
n}$ be the set of all $m\times n$ real matrices, $\dbS^n$ be the set
of all $n\times n$ symmetric real matrices, and for any
$M=(m_{ij})\in\dbR^{m\times n}$, $M^T$ stands for its transpose and
let
$$|M|=\(\sum_{i,j}m_{ij}^2\)^{1\over2},\qq\forall M\in\dbR^{m\times n}.$$
For $X=\dbR^n,\dbR^{n\times m}$, etc., let
$$\ba{ll}
\ns\ds L^2(a,b;X)=\Big\{f:[a,b]\to X\bigm|f(\cd)\hb{ is measurable,
}\int_a^b|f(t)|^2dt<+\infty\Big\},\\
\ns\ds L^\infty(a,b;X)=\Big\{f:[a,b]\to X\bigm|f(\cd)\hb{ is measurable, }\esssup_{t \in[a,b]}|f(t)|<\infty\Big\},\\
\ns\ds C([a,b];X)=\Big\{f:[a,b]\to X\bigm|f(\cd)\hb{ is continuous}\Big\},\ea$$
and
$$\ba{ll}
\ns\ds L^p_{\cF_t}(\Om;X)=\Big\{\xi:\Om\to X\bigm|\xi\hb{ is $\cF_t$-measurable, }\dbE|\xi|^p<\infty\Big\},\qq p\ge1,\\
\ns\ds L^\infty_{\cF_t}(\Om;X)=\Big\{\xi\in L^1_{\cF_t}(\Om;X)\bigm|\esssup_{\o\in\Om}|\xi(\o)|<\infty\Big\},\\
\ns\ds L_\dbF^p(a,b;X)=\Big\{f:[a,b]\times\Omega\to
X\bigm|f(\cd)\hb{ is
$\dbF$-adapted, }\dbE\int_a^b|f(t)|^pdt<+\infty\Big\},\q p\ge1,\\
\ns\ds L_\dbF^\infty(a,b;X)=\Big\{f(\cd)\in
L^1_\dbF(a,b;X)\bigm|\esssup_{(t,
\omega)\in[a,b]\times\Omega}|f(t,\o)|<\infty\Big\},\\
\ns\ds C_\dbF\big([a,b];L^p_{\cF_T}(\Om;X)\big)=\Big\{f(\cd) \in
L^p_\dbF(a,b:X)\bigm|t\mapsto f(t)\hb{ is continuous, }\sup_{t\in[a,b]}\dbE|f(t)|^p<\infty\big\},\\
\ns\ds L^p_\dbF(\Om;C([a,b];X))=\Big\{f(\cd)\in
L^p_\dbF(a,b;X)\bigm|\dbE\[\sup_{t\in[a,b]}|f(t)|^p\]<+\infty\Big\},\qq
p\ge1.\ea$$
In the definition of $C_\dbF([a,b];L^p_{\cF_T}(\Om;X))$, $t\mapsto
f(t)$ is continuous means that
$$\lim_{t\to t_0}\dbE|f(t)-f(t_0)|^p=0,\qq\forall t_0\in[a, b].$$
Now, let us introduce the following sets:
$$\ba{ll}
\ns\ds\cU[0,T]=\Big\{u:[0,T]\times\Om \to U\bigm|u(\cd)\in
L^2_\dbF(0,T;\dbR^m)\Big\},\\
\ns\ds\cT[0,T]=\Big\{\t:\Om\to[0,T]\bigm|\t\hb{ is an
$\dbF$-stopping
time}\Big\}.\ea$$
%
%
Hereafter, $U\subseteq\dbR^m$ is assumed to be convex and closed.
Any $u(\cd)\in\cU[0,T]$ is called a {\it regular admissible control}
while $\t\in\cT[0,T]$ is called an {\it admissible stopping time}.
Under some mild conditions, for any $x_1\in\dbR^{n_1}$, and
$(u(\cd),\t)\in\cU[0,T]\times\cT[0,T]$, (\ref{state11}) admits a
unique strong solution $X(\cd)\equiv X(\cd\,;x_1,u(\cd),\t)$, and the
cost functional (\ref{cost}) is well-defined. Having this, we can
pose the following problem.

\ms

\fbox{\textbf{Problem (MLQ)}}\ \  For given $x_1\in\dbR^{n_1}$, find
a $(\bar u(\cd),\bar\t)\in\cU[0,T]\times\cT[0,T]$ such that
$$
J(x_1;\bar
u(\cd),\bar\t)=\inf_{(u(\cd),\t)\in\cU[0,T]\times\cT[0,T]}J(x_1;u(\cd),\t)
\equiv V(x_1).
$$
Any $(\bar u(\cd),\bar\t)\in\cU[0,T]\times\cT[0,T]$ satisfying the
above is called an optimal control pair, $\bar X\1n(\cd)\1n\equiv\2n
X\1n(\cd\,;x_1,\bar u(\cd),\bar\t)$ is called the corresponding
optimal trajectory, $(\bar X(\cd),\bar u(\cd),\bar\t)$ is called an
{\it optimal triple}. In the above, MLQ problem stands for {\it
mixed linear-quadratic} problem, in which, one has a usual control
$u(\cd)$ mixed with a control $\t$ of stopping time. We have the
following points to the above MLQ problem formulation.

\begin{itemize}

\item A special feature of Problem (MLQ) is that in minimizing the
cost functional, one needs to select a regular control $u(\cd)$ from
$\cU[0,T]$ and at the same time, one has to find the best time
$\bar\t$ to initiate or trigger the whole system. Note that the initial
value $X(\t)$ of $X(\cd)$ at $\t$ depends on the value of
$X_1(\t-0)$, which in turn depends on the regular control on
$[0,\t)$. From this viewpoint, the Problem (MLQ) is some kind of
combination of a usual stochastic optimal control and an optimal
stopping time problems. Similar problems have been investigated in
literature, including {\O}ksendal and Sulem \cite{BS} where some
optimal resource extraction optimization problem was addressed. By
applying the dynamic programming method, the optimal policy can be
characterized by some Hamilton-Jacobi-Bellman variational
inequalities, see Krylov \cite{K}, for relevant treatment. In
contrast, here, we aim to investigate the problem by the variational
method which could lead to Pontryagin type maximum principle.

\item Consider the following simple but illustrating example, from which
we can see the significant difference between Problem (MLQ) and
other relevant ones when applying the possible perturbation method.
Suppose the controlled state equation is given by
$$\left\{\ba{ll}
\ns\ds dX_1(t)=a_1X_1(t)dW(t),\qq t\in[0,T],\\
\ns\ds dX_2(t)=a_2X_2(t)dW(t),\qq t\in[\t,T],\\
\ns\ds X_2(\t)=KX_1(\t),\ea\right.$$
where $a_1, a_2, K$ are some constants. Let $\bar X(\cd)\equiv(\bar
X_1(\cd),\bar X_2(\cd))$ be the solution corresponding to $\bar\t$.
Introduce a perturbation on $\bar\t$ of the form:
$\t^\rho=\bar\t+\rho\t$, $\rho>0$, with $\t$ being another stopping
time. Let the solution corresponding to $\t^\rho$ be
$X^\rho(\cd)\equiv(X^\rho_1(\cd),X^\rho_2(\cd))$. Then
$X^\rho_1(\cd)=\bar X_1(\cd)$ is independent of $\rho$, and for
$t\in[\t^\rho,T]$,
$$X_2^\rho(t)-\bar X_2(t)=K\(\bar X_1(\t^\rho)-\bar X_1(\bar\t)\)+\int_{\t^\rho}^ta_2 \left(X_2^\rho(s)-\bar X_2(s)\right)dW(s)-\int_{\bar\t}^{\t^\rho}a_2\bar X_2(s)dW(s).$$
The presence of the term $\int_{\bar\t}^{\t^\rho}a_2\bar
X_2(s)dW(s)$ makes the first-order Taylor expansion in convex
variation failed to work. This is mainly due to the fact that
$$\dbE\left(\int_{\bar\t}^{\t^\rho}a_2\bar X_2(s)dW(s)\right)^2=\dbE\int_{\bar\t}^{\t^\rho}(a_2\bar X_2(s))^{2}ds,$$
which, in rough sense, suggests $\int_{\bar\t}^{\t^\rho}a_2\bar
X_2(s)dW(s)$ be of order $\sqrt\rho$ instead $\rho$. On the other
hand, it is not feasible to apply the second-order Taylor expansion to
introduce the second-order variational equation (as suggested by
Peng \cite{P}, Yong and Zhou \cite{YZ} etc.). This is mainly because
the time horizon on which $X_2(\cd)$ is defined depends on the
selection of $\t$, thus the spike variation method cannot be applied
here either.

\item Problem (MLQ) also differs from the well-studied stochastic impulse
control, since as the time passes $\t$, instead of having a jump
for the state as a usual impulse control does, our controlled system
changes the dimension of the state (from $X_1(\cd)$ to $X(\cd)$).

\end{itemize}

In summary, the involvement of $\t$ into the control variable makes
the Problem (MLQ) essentially different from other classical optimal
control problems, and the standard perturbation jointly on $(\bar
u(\cd), \bar\t)$ is not workable directly. Keep this in mind, in
this paper, we take the following strategy to study Problem (MLQ):
we first connect the Problem (MLQ) into some random-duration linear
quadratic (RLQ, for short) optimal control problem, and an optimal
time (OT, for short) problem to the associated Riccati equations. By
Problem (RLQ), we can obtain some necessary condition for the
regular optimal control $\bar u(\cd)$; by Problem (OT), we can
obtain some necessary condition satisfied by the optimal time
$\bar\t$. Next, by solving Problem (RLQ) and Problem (OT)
consecutively, we can solve
the original Problem (MLQ).\\

\ms

The rest of this paper is organized as follows. In Section 2, we get
a stochastic maximum principle for a little more general two-stage
random-duration optimal control problems. Based on it, Section 3 is
devoted to a study of the random-duration linear quadratic optimal
control problems. The state feedback optimal control is derived via
some stochastic Riccati-type equations and the optimal cost
functional is also calculated explcitly. In Section 4, an equivalence between
Problem (MLQ) and Problems (RLQ)--(OT) is established. In Section 5,
for the case of one-dimension with constant coefficients, we
characterize the optimal time $\bar\t$.

\section{Random-Duration Optimal Control Problem}

In this section, we consider the following controlled stochastic
differential equation (SDE, for short)
\bel{2.2'}\left\{\ba{ll}
\ns\ds dX_1(t)=b^1(t,X_1(t),u(t))dt+\si^1(t,X_1(t),u(t))dW(t),\q t\in[0,\t),\\
\ns\ds dX(t)=b(t,X(t),u(t))dt+\si(t,X(t),u(t))dW(t),\q t\in[\t,T],\\
\ns\ds X_1(0)=x_1, \qq X(\t)=\Phi(\t,X_1(\t)),\ea\right.\ee
where $\t\in\cT(0,T]$ is some fixed stopping time and
$u(\cd)\in\cU[0,T]$ is an admissible control. The cost functional is
\bel{functional}J^\t(x_1;u(\cd))=\dbE\[\int_0^\t
g^1(t,X_1(t),u(t))dt+\int_\t^Tg(t,X(t),u(t))dt+h^1(\t,X_1(\t))+h(X(T))\].\ee
Consider the following random-duration optimal control
(ROC, for short) problem:\\

\fbox{\textbf{Problem (ROC)}} For $x_1\in\dbR^{n_1}$ and
$\t\in\cT[0,T]$, find a $\bar u(\cd)\in\cU[0,T]$ such that
\bel{3.1}J^\t(x_1;\bar
u(\cd))=\inf_{u(\cd)\in\cU[0,T]}J^\t(x_1;u(\cd))\equiv V^\t(x_1).\ee

The following basic assumptions will be in force:\\

\textbf{(H2.1)} Let
$$\ba{ll}
\ns\ds b,\si:[0,T]\times\dbR^n\times U\times\Om\to\dbR^n,
\q b^1,\si^1:[0,T]\times\dbR^{n_1}\times U\times\Om\to\dbR^{n_1},\\
\ns\ds\Phi:[0,T]\times\dbR^{n_1}\times\Om\to\dbR^n\ea$$
be measurable such that
$$\ba{ll}
\ns\ds(t,\o)\mapsto(b(t,x,u,\o),\si(t,x,u,\o)),\\
\ns\ds(t,\o)\mapsto(b^1(t,x,u,\o),\si^1(t,x,u,\o),\F(t,x_1,\o)),\ea$$
are progressively measurable,
$$\ba{ll}
\ns\ds(x,u)\mapsto(b(t,x,u,\o),\si(t,x,u,\o)),\\
\ns\ds(x_1,u)\mapsto(b^1(t,x_1,u,\o),\si^1(t,x_1,u,\o),\Phi(t,x_1,\o)),\ea$$
are continuously differentiable, and for some constant $L>0$,
$$\ba{ll}
\ns\ds|b_x(t,x,u)|+|b_u(t,x,u)|+|b^1_{x_1}(t,x_1,u)|+|\si^1_{x_1}(t,x_1,u)| +|\Phi_{x_1}(t,x_1)|\le L,\\
\ns\ds\qq\qq\qq\qq\qq\qq\qq\qq\qq\forall(t,x,u)\in[0,T]\times\dbR^n\times
U,~\as,\ea$$
and
$$
|b(t,0,u)|+|\si(t,0,u)|+|b^1(t,0,u)|+|\si^1(t,0,u)|+|\Phi(t,0)|\le
L,\qq(t,u)\in[0,T]\times U,~\as
$$
\\

\textbf{(H2.2)} Let
$$\ba{ll}
\ns\ds g:[0,T]\times\dbR^n\times U\times\Om\to\dbR,
\q g^1:[0,T]\times\dbR^{n_1}\times U\times\Om\to\dbR,\\
\ns\ds h:\dbR^n\times\Om\to\dbR,\q
h^1:[0,T]\times\dbR^{n_1}\times\Om\to\dbR\ea$$
be measurable such that
$$\ba{ll}
\ns\ds(t,\o)\mapsto g(t,x,u,\o)),\qq(t,\o)\mapsto g^1(t,x,u,\o)\ea$$
are progressively measurable,
$$\ba{ll}
\ns\ds(x,u)\mapsto g(t,x,u,\o),\qq(x_1,u)\mapsto
g^1(t,x_1,u,\o),\ea$$
$$(x_1,u)\mapsto(b^1(t,x_1,u,\o),\si^1(t,x_1,u,\o),\Phi(t,x_1,\o)),$$
are continuously differentiable, and for some constant $L>0$,
$$\ba{ll}
\ns\ds|g_x(t,x,u)|+|g_u(t,x,u)|+|h_x(x)|\le L(1+|x|+|u|),\\
\ns\ds|g^1_{x_1}(t,x_1,u)|+|g^1_u(t,x_1,u)|+|h^1_{x_1}(t,x_1)|\le L(1+|x_1|+|u|),\\
\ns\ds\qq\qq\qq\qq\qq\qq\qq\forall(t,x,u)\in[0,t]\times\dbR^n\times
U,\ea$$
and
$$|g(t,0,u)|+|g^1(t,0,u)|+|h^1(t,0)|\le L,\qq\forall(t,u)\in[0,T]\times U.$$

By some standard arguments, we see that under assumptions
{\bf(H2.1)}--{\bf(H2.2)}, for any $x_1\in\dbR^{n_1}$, and
$(u(\cd),\t)\in\cU[0,T]\times\cT[0,T]$, the state equation
$(\ref{2.2'})$ admits a unique solution $X(\cd\,;x_1,u(\cd),\t)$ and
the cost functional $J^\t(x_1;u(\cd))$ is well-defined. Therefore,
Problem (ROC) makes sense. Suppose $(\bar X(\cd),\bar u(\cd))$ is an
optimal pair of Problem (ROC), depending on $(x_1,\t)$. For any
$v(\cd)\in\cU[0,T]$, denote
$$u(\cd)=v(\cd)-\bar u(\cd).$$
By the convexity of $U$, we know that
$$u^\rho(\cd)=\bar u(\cd)+\rho u(\cd)=(1-\rho)\bar u(\cd)+\rho v(\cd)\in\cU[0,T],\qq\forall\rho\in(0,1).$$
By the optimality of $(\bar X(\cd),\bar u(\cd))$, we have
$$\liminf_{\rho\to0}{J^\t(x_1;u^\rho(\cd))-J^\t(x_1;\bar u(\cd))\over\rho}\ge0.$$
Making use of some similar arguments in \cite{Mou-Yong}, we have the
following result.

\ms

\bf Lemma 2.1. \sl Suppose
\emph{\textbf{(H2.1)}}--\emph{\textbf{(H2.2)}} hold. Let $(\bar
X(\cd),\bar u(\cd))$ be an optimal pair of Problem {\rm(ROC)}. Then
\bel{variational-inequality}\ba{ll}
\ns\ds\dbE\Big\{\int_0^\t\big[g^1_{x_1}(t,\bar X_1(t),\bar u(t))\xi_1(t)+g_u^1(t,\bar X_1(t),\bar u(t))u(t)\big]dt+h^1_{x_1}(\t,\bar X_1(\t))\xi_1(\t)\\
\ns\ds\qq+\int_\t^T\big[g_x(t,\bar X(t),\bar u(t))\xi(t) +g_u(t,\bar
X(t),\bar u(t))u(t)\big]dt+h_x(\bar X(T))\xi(T)\Big\}\ge0,\ea\ee
where $\xi(\cd)\equiv(\xi_1(\cd),\xi_2(\cd))$ is the solution to the
following variational system:
$$\left\{\ba{ll}
\ns\ds d\xi_1(t)=\big[b^1_{x_1}(t,\bar X_1(t),\bar u(t))\xi_1(t)+b^1_u(t,\bar X_1(t),\bar u(t))u(t)\big]dt \\
\ns\ds\qq\qq+\big[\si^1_{x_1}(t,\bar X_1(t),\bar u(t))\xi_1(t)+\si^1_u(t,\bar X_1(t),\bar u(t))u(t)\big]dW(t),\q t\in[0,\t), \\
\ns\ds d\xi(t)=\big[b_x(t,\bar X(t),\bar u(t))\xi(t)+b_u(t,\bar X(t),\bar u(t))u(t)\big]dt \\
\ns\ds\qq\qq+\big[\si_x(t,\bar X(t),\bar u(t))\xi(t)+\si_u(t,\bar X(t),\bar u(t))u(t)\big]dW(t),\q t\in[\t,T],\\
\ns\ds\xi_1(0)=0,\qq\xi(\t)=\Phi_{x_1}(\t,\bar X_1(\t))\xi_1(\t).
\ea\right.$$

\ms

\rm

The following is a Pontryagin type maximum principle.

\ms

\bf Theorem 2.2. \sl Suppose
\emph{\textbf{(H2.1)}}--\emph{\textbf{(H2.2)}} hold. Let $(\bar
X(\cd),\bar u(\cd))$ be an optimal pair of Problem {\rm (ROC)}. Then
the following two backward stochastic differential equations (BSDEs,
for short) admit unique adapted solutions $(p(\cd),q(\cd))$ and
$(p_1(\cd),q_1(\cd))$:
\bel{Hamiltonian-2} \left\{\ba{ll}
\ns\ds dp(t)=-\big[b_x(t,\bar X(t),\bar u(t))^Tp(t)+\si_x(t,\bar X(t),\bar u(t))^Tq(t)\\
\ns\ds\qq\qq\q-g_x(t,\bar X(t),\bar u(t))^T\big]dt+q(t)dW(t),\qq t\in[\t,T],\\
\ns\ds p(T)=-h_x(\bar X(T))^T,\ea\right.\ee
\bel{Hamiltonian-1}\left\{\1n\ba{ll}
\ns\ds dp_1(t)=-\big[b^1_{x_1}(t,\bar X_1(t),\bar u(t))^Tp_1(t)+\si^1_{x_1}(t,\bar X_1(t),\bar u(t))^Tq_1(t)\\
\ns\ds\qq\qq\q-g^1_{x_1}(t,\bar X_1(t),\bar u(t))^T\big]dt+q_1(t)dW(t),\qq t\in[0,\t],\\
\ns\ds p_1(\t)=-h^1_{x_1}(\t,\bar X_1(\t))^T+\Phi_{x_1}(\t,\bar
X_1(\t))^Tp(\t).\ea\right.\ee
Moreover, the following variational inequalities hold
\bel{optimal-u}\left\{\ba{ll}
\ns\ds\[p_1(t)^Tb^1_u(t,\bar X_1(t),\bar u(t))+q_1(t)^T\si_u^1(t,\bar X_1(t),\bar u(t))-g^1_u(t,\bar X_1(t),\bar u(t))\][v-\bar u(t)]\le0,\\
\ns\ds\qq\qq\qq\qq\qq\qq\qq t\in[0,\t),\q v\in U,\q\as,\\
\ns\ds\[p(t)^Tb_u(t,\bar X(t),\bar u(t))+q(t)^T\si_u(t,\bar X(t),\bar u(t))-g_u(t,\bar X(t),\bar u(t))\][v-\bar u(t)]\le0,\\
\ns\ds\qq\qq\qq\qq\qq\qq\qq t\in[\t,T],\q v\in U,\q\as\ea\right.\ee

\ms

\rm

\it Proof. \rm Applying It\^o formula to $\lan
p_1(\cd),\xi_1(\cd)\ran$ and $\lan p(\cd),\xi(\cd)\ran$,
respectively, we have
$$\ba{ll}
\dsp\dbE\[\lan-h^1_{x_1}(\t,\bar X_1(\t))^T+\Phi_{x_1}(\t,\bar
X_1(\t))^T
p(\t),\xi_1(\t)\ran\]=\dbE\lan p_1(\t),\xi_1(\t)\ran\\
\ns\ds=\dbE\int_0^\t\[\lan-\big[b^1_{x_1} (t,\bar X_1,\bar
u)^Tp_1+\si^1_{x_1}(t,\bar X_1,\bar u)^Tq_1
-g^1_{x_1}(t,\bar X_1,\bar u)^T\big],\xi_1\ran\\
\dsp\qq\qq+\lan p_1,b^1_{x_1}(t,\bar X_1,\bar u)\xi_1+b^1_u(t,\bar X_1,\bar u)u\ran+\lan q_1,\si^1_{x_1}(t,\bar X_1,\bar u)\xi_1+\si^1_u(t,\bar X_1,\bar u)u\ran\]dt\\
\dsp=\dbE\int_0^\t\[g^1_{x_1}(t,\bar X_1,\bar u)\xi_1+\lan
b^1_u(t,\bar X_1,\bar u)^Tp_1+\si^1_u(t,\bar X_1,\bar
u)^Tq_1,u\ran\]dt,\ea$$
and
$$\ba{ll}
\dsp\dbE\lan-h_x(\bar X(T))^T,\xi(T)\ran=\dbE\lan p(T),\xi(T)\ran\\
\ns\ds=\dbE\Big\{\lan p(\t),\xi(\t)\ran+\int_\t^T\[\lan-\big[b_x(t,\bar X,\bar u)^Tp+\si_x(t,\bar X,\bar u)^Tq -g_x(t,\bar X,\bar u)^T\big],\xi\ran\\
\dsp\qq\qq\qq\qq\qq+\lan p,b_x(t,\bar X,\bar u)\xi+b_u(t,\bar X,\bar u)u\ran+\lan q,\si_x(t,\bar X,\bar u)\xi+\si_u(t,\bar X,\bar u)u\ran\]dt\Big\}\\
\dsp=\dbE\Big\{\lan p(\t),\Phi_{x_1}(\t,\bar
X_1(\t))\xi_1(\t)\ran+\int_\t^T\[g_x(t,\bar X,\bar u)\xi+\lan
b_u(t,\bar X,\bar u)^Tp+\si_u(t,\bar X,\bar
u)^Tq,u\ran\]dt\Big\}.\ea$$
Then, we obtain
$$\ba{ll}
\ns\ds0\le\dbE\Big\{\int_0^\t\big[g^1_{x_1}(t,\bar X_1,\bar u)\xi_1+g_u^1(t,\bar X_1,\bar u)u\big]dt+h^1_{x_1}(\t,\bar X_1(\t))\xi_1(\t)\\
\ns\ds\qq\qq+\int_\t^T\big[g_x(t,\bar X,\bar u)\xi+g_u(t,\bar X,\bar
u)u\big]dt+h_x(\bar X(T))\xi(T)\Big\}\\
\ns\ds\q=\dbE\Big\{\int_0^\t\lan-b^1_u(t,\bar X_1,\bar u)^Tp_1-\si^1_u(t,\bar X_1,\bar u)^Tq_1+g^1_u(t,\bar X_1,\bar u)^T,u\ran dt+\lan\Phi_{x_1}(\t,\bar X_1(\t))^Tp(\t),\xi_1(\t)\ran \\
\ns\ds\qq\q+\int_\t^T\lan-b_u(t,\bar X,\bar u)^Tp-\si_u(t,\bar X,\bar u)^Tq+g_u(t,\bar X,\bar u)^T,u\ran dt-\lan p(\t),\Phi_{x_1}(\t,\bar X_1(\t))\xi_1(\t)\ran\Big\}\ea$$
$$\ba{ll}
\ns\ds\q=\dbE\Big\{\int_0^\t\lan-b^1_u(t,\bar X_1,\bar u)^Tp_1-\si^1_u(t,\bar X_1,\bar u)^Tq_1+g^1_u(t,\bar X_1,\bar u)^T,v-\bar u\ran dt \\
\ns\ds\qq\q+\int_\t^T\lan-b_u(t,\bar X,\bar u)^Tp-\si_u(t,\bar
X,\bar u)^Tq+g_u(t,\bar X,\bar u)^T,v-\bar u\ran
dt\Big\}.~~~~~~~~~~~~~~~~~~~~~~~~~~~~~~~~~~\ea$$
Note that $v(\cd)\in\cU[0,T]$ is arbitrary, we therefore have
$(\ref{optimal-u})$. \endpf

\section{Random-Duration Linear Quadratic Problem}
For any $(x_1,\t)\in\dbR^{n_1}\times\cT[0,T]$, we consider the
following controlled linear system:
\bel{3.0}\left\{\ba{ll}
\dsp dX_1(t)=\big[A_1(t)X_1(t)+B_1(t)u(t)\big]dt+\big[C_1(t)X_1(t)+D_1(t)u(t)\big]dW(t),\q t\in[0,\t),\\[1mm]
\dsp dX(t)=\big[A(t)X(t)+B(t)u(t)\big]dt+\big[C(t)X(t)+D(t)u(t)\big]dW(t),\q t\in[\t,T], \\ [1mm]
\dsp X_1(0)=x_1,\qq X(\t)=K(\t)X_1(\t-0),\ea\right.\ee
with quadratic cost functional as follows:
\bel{3.0'}\ba{ll}
\dsp J^\t(x_1;u(\cd))={1\over2}\dbE\Big\{\int_0^\t\big[\lan Q_1(t)X_1(t),X_1(t)\ran+\lan R_1(t)u(t),u(t)\ran\big]dt+\lan G_1(\t)X_1(\t),X_1(\t)\ran \\
\ns\ds\qq\qq\qq\qq\qq+\int_\t^T\big[\lan Q(t)X(t),X(t)\ran+\lan R(t)u(t),u(t)\ran\big]dt+\lan GX(T),X(T)\ran\Big\}.
\ea\end{equation}
Now we pose the following problem.

\ms

\fbox{\textbf{Problem (RLQ)}} For given
$(x_1,\t)\in\dbR^{n_1}\times\cT[0,T]$, find a $\bar
u(\cd)\in\cU[0,T]$ such that
\begin{equation}\label{3.1}
J^\t(x_1;\bar u(\cd))=\inf_{u(\cd)\in\cU[0,T]}J^\t(x_1;u(\cd))\equiv
V^\t(x_1).
\end{equation}

We call the above a \emph{random-duration linear quadratic} (RLQ,
for short) problem. For the above problem, we introduce the
following hypothesis.

\ms

\textbf{(H3.1)} The following holds:
$$\left\{\ba{ll}
\ns\ds A_1(\cd),C_1(\cd)\in L^\infty_\dbF(0,T;\dbR^{n_1\times n_1}),\qq A(\cd),C(\cd)\in L^\infty_\dbF(0,T;\dbR^{n\times n}),\\
\ns\ds B_1(\cd),D_1(\cd)\in L^\infty_\dbF(0,T;\dbR^{n_1\times m}),\qq B(\cd),D(\cd)\in L^\infty_\dbF(0,T;\dbR^{n\times m}), \\
\ns\ds Q_1(\cd)\in L^\infty_\dbF(0,T;\dbS^{n_1}),\q Q(\cd)\in
L^\infty_\dbF(0,T;\dbS^n),\q R_1(\cd), R(\cd)\in L^\infty_\dbF(0,T;\dbS^m),\\
\ns\ds G_1(\cd)\in C_\dbF([0,T];L^\infty_{\cF_T}(\Om;\dbS^{n_1})),\q
G\in L^\infty_\dbF(\Om;\dbS^n),\q K(\cd)\in
C_\dbF([0,T];L^\infty_{\cF_T}(\Om;\dbR^{n\times n_1})).\ea\right.$$
Moreover, for some $\d>0$,
\bel{}\left\{\ba{ll}
\ns\ds R_1(t)\ge\d I_{n_1},\q R(t)\ge\d I_n,\\
\ns\ds Q_1(t),G_1(t)\ge0,\q Q(t)\ge0,\q G\ge0,\ea\right.\qq t\in[0,T],~\as\ee
It is clear that under {\bf(H3.1)}, for any given
$(x_1,\t)\in\dbR^{n_1}\times\cT[0,T]$, the map $u(\cd)\mapsto
J(x_1;u(\cd),\t)$ is convex and coercive. Therefore, Problem (RLQ)
admits a unique optimal control $\bar u(\cd)\in\cU[0,T]$. Now, let
$(\bar X(\cd),\bar u(\cd))$ be the optimal pair of Problem (RLQ),
depending on $(x_1,\t)\in\dbR^{n_1}\times\cT[0,T]$. By Theorem
$2.2$, on $[\t,T]$ the optimal pair $(\bar X(\cd),\bar
u(\cd))$ satisfies the following:
\bel{3.2}\left\{\ba{ll}
\ns\ds d\bar X(t)=\big[A(t)\bar X(t)+B(t)\bar u(t)\big]dt+\big[C(t)\bar X(t)+D(t)\bar u(t)\big]dW(t),\\
\ns\ds dp(t)=-\big[A(t)^Tp(t)+C(t)^Tq(t)-Q(t)\bar X(t)\big]dt+q(t)dW(t),\\
\ns\ds p(T)=-G\bar X(T),\qq \bar X(\t)=K(\t)\bar X_1(\t),\\
\ns\ds B(t)^Tp(t)+D(t)^Tq(t)-R(t)\bar u(t)=0.\ea\right.\ee
The above is a coupled FBSDE in random duration. To solve it, we let
$$p(t)=-P(t)\bar X(t),\qq t\in[0,T],$$
for some $P(\cd)$ satisfying
$$\left\{\ba{ll}
\ns\ds dP(t)=\G(t)dt+\L(t)dW(t),\qq t\in[0,T],\\
\ns\ds P(T)=G.\ea\right.$$
Then apply It\^o's formula, we have
$$\ba{ll}
\ns\ds\big(-A^TP\bar X+C^Tq-Q\bar X\big)dt-qdW=-dp=d\big(P\bar X\big)\\
\ns\ds=\Big\{\G\bar X+P\big(A\bar X+B\bar u\big)+\L\big(C\bar X+D\bar u\big)\Big\}dt+\Big\{\L\bar X+P\big(C\bar X+D\bar u\big)\Big\}dW\\
\ns\ds=\Big\{(\G+PA+\L C)\bar X+(PB+\L D)\bar u\Big\}dt+\Big\{(\L+PC)\bar X+PD\bar u\Big\}dW.\ea$$
Hence,
$$q=-(\L+PC)\bar X-PD\bar u.$$
Then
$$R\bar u=B^Tp+D^Tq=-B^TP\bar X-D^T(\L+PC)\bar X-D^TPD\bar u.$$
This implies
$$\bar u=-(R+D^TPD)^{-1}\big(B^TP+D^T\L+D^TPC\big)\bar X.$$
Substituting $\bar u$ into the expression of the above $q$ yields
$$q=-(\L+PC)\bar X+PD(R+D^TPD)^{-1}(B^TP+D^T\L+D^TPC)\bar X.$$
Moreover,
$$\ba{ll}
\ns\ds0=(A^TP+Q)\bar X-C^Tq+(\G+PA+\L C)\bar X+(PB+\L D)\bar u\\
\ns\ds\q=(\G+PA+A^TP+\L C+Q)\bar X \\
\ns\ds\qq-C^T\[-(\L+PC)\bar X+PD(R+D^TPD)^{-1}(B^TP+D^T\L+D^TPC)\bar X\] \\
\ns\ds\qq-(PB+\L D)(R+D^TPD)^{-1}\big(B^TP+D^T\L+D^TPC\big)\bar X \\
\ns\ds\q=\(\G+PA+A^TP+C^TPC+\L C+C^T\L+Q \\
\ns\ds\qq-(PB+\L D+C^TPD)(R+D^TPD)^{-1}\big(B^TP+D^T\L+D^TPC\big)\)\bar X.\ea$$
Therefore, we take
\bel{Gamma}\ba{ll}
\ns\ds\G=-\[PA+A^TP+C^TPC+\L C+C^T\L+Q\\
\ns\ds\qq\qq-(PB+\L D+C^TPD)(R+D^TPD)^{-1}\big(B^TP+D^T\L+D^TPC\big)\].\ea\ee
Consequently, the corresponding Riccati equation reads
\bel{P-Riccati0}\left\{\ba{ll}
\ns\ds dP=-\[PA+A^TP+C^TPC+\L C+C^T\L+Q\\
\ns\ds\qq\qq-(PB+\L D+C^TPD)(R+D^TPD)^{-1}\big(B^TP+D^T\L+D^TPC\big)\]dt\\
\ns\ds\qq\qq+\L dW,\qq t\in[0,T],\\
\ns\ds P(T)=G,\qq R+D^TPD>0.\ea\right.\ee
On $[0,\t]$ the optimal pair $(\bar X_1(\cd),\bar u_1(\cd))$ satisfies
\bel{3.3FBSDE}\left\{\3n\ba{ll}
\ns\ds d\bar X_1\1n=\1n\big[A_1\bar X_1\1n+\1n B_1R_1^{-1}B_1^Tp_1\1n+\1n B_1R_1^{-1}D_1^Tq_1\big]dt\1n+\1n\big[C_1\bar X_1\1n+\1n D_1R_1^{-1}B_1^Tp_1\1n+\1n D_1R_1^{-1}D_1^Tq_1\big]dW(t),\\
\ns\ds dp_1=-\big[A_1^Tp_1+C_1^Tq_1-Q_1\bar X_1\big]dt+q_1dW(t),\\
\ns\ds\bar X_1(0)=x_1,\qq p_1(\t)=-\[K(\t)^TP(\t)K(\t)+G_1(\t)\]\bar X_1(\t), \\
\ns\ds R_1\bar u_1=B_1^Tp_1+D_1^Tq_1.\ea\right.\ee
Similar to the above, the corresponding Riccati equation is
\bel{P_1-Riccati0}\left\{\ba{ll}
\ns\ds dP_1=-\[P_1A_1+A_1^TP_1+C_1^TP_1C_1+\L_1C_1+C_1^T\L_1+Q_1\\
\ns\ds\qq\qq-(P_1B_1+\L_1D_1+C_1^TP_1D_1)(R_1+D_1^TP_1D_1)^{-1}\big(B_1^TP_1+D_1^T\L_1+D_1^TP_1C_1\big)\]dt\\
\ns\ds\qq\qq+\L_1dW,\qq t\in[0,\t],\\
\ns\ds P_1(\t)=K(\t)^TP(\t)K(\t)+G_1(\t),\qq R_1+D_1^TP_1D_1>0.\ea\right.\ee

Following \cite{T}, we know that under \textbf{(H3.1)}, Riccati
equations (\ref{P-Riccati0}) and (\ref{P_1-Riccati0}) admit unique
adapted solutions on $[0,T]$ and $[0,\t]$ respectively. The
following is a verification theorem.

\ms

\bf Theorem 3.1. \sl Let \textbf{(H3.1)} hold. Let $P(\cd)$ and
$P_1(\cd)$ be the solutions of Riccati equations
$(\ref{P-Riccati0})$ and $(\ref{P_1-Riccati0})$ respectively. Then
Problem (RLQ) has an optimal control with state feedback form as
follows:
\bel{bar u}\bar u(t)=\begin{cases}
-\Psi(t)\bar X(t),& t \in [\t,T],\\
-\Psi_1(t)\bar X_1(t), & t \in [0,\t),
\end{cases}\ee
where
\bel{Psi}\left\{\3n\ba{ll}
\ns\ds\Psi(t)=\[R(t)+D(t)^TP(t)D(t)\]^{-1}\[B(t)^TP(t)+D(t)^TP(t)C(t)+D(t)^T\L(t)\],\\
\ns\ds\Psi_1(t)=\[R_1(t)+D_1(t)^TP_1(t)D_1(t)\]^{-1}\[B_1(t)^TP_1(t)+D_1(t)^TP_1(t)C_1(t)+D_1(t)^T\L_1(t)\],\ea\right.\ee
Moreover, the optimal value of the cost functional is given by
\bel{J(bar u)}J^\t(x_1;\bar u(\cd))={1\over2}\lan
P_1(0)x_1,x_1\ran=V^\t(x_1).\ee

\ms

\it Proof. \rm For any  $u(\cd)\in\cU[0,T]$, let $X(\cd)$ be the
corresponding state process. Applying It\^o formula to $\lan
P(\cd)X(\cd),X(\cd)\ran$ on the interval $[\t,T]$, we obtain (let
$\G$ be defined by (\ref{Gamma}))
$$\ba{ll}
\dsp\dbE\big[\lan GX(T),X(T)\ran-\lan P(\t)X(\t),X(\t)\ran\big]
\\ [2mm]
\ns\ds=\dbE\int_\t^T\(\lan(\G+PA+A^TP+C^TPC+\L C+C^T\L)X,X\ran \\
\ns\ds\qq\qq\qq+2\lan(B^TP+D^T\L+D^TPC)X,u\ran+\lan D^TPDu,u\ran\)dt\\
\dsp=\dbE\Big\{\int_\t^T\big[2\lan(B^TP+D^TPC+D^T\L)X,u\ran -\lan QX,X\ran+\lan D^TPDu,u\ran \\ [2mm]
\dsp\qq\qq+\lan(PB+C^TPD+\L D)(R+D^TPD)^{-1}(B^TP+D^TPC+D^T\L)X,X\ran\big]dt\Big\}. \ea$$
It follows that:
$$\ba{ll}
\dsp J^{\t}(x_1;u(\cd))-\frac{1}{2}\dbE \lan
P(\t)X(\t),X(\t)\ran\\ [2mm]
\dsp={1\over2}\dbE\Big\{\int_0^\t\(\lan Q_1X_1,X_1\ran+\lan R_1u,u\ran\)dt+\lan G_1(\t)X_1(\t),X_1(\t)\ran \\
\ns\ds\qq+\int_\t^T\(\lan QX,X\ran+\lan Ru,u\ran\)dt+\lan GX(T),X(T)\ran\\
\ns\ds\qq-\lan GX(T),X(T)\ran\1n+\1n\int_\t^T\big[2\lan
(B^T\1n P\1n+\1n D^T\1n PC\1n+\1n D^T\1n\L)X,u\ran\1n-\1n\lan QX,X\ran\1n+\1n\lan D^T\1n PDu,u\ran \\
[2mm]
\dsp\qq+\lan(PB+C^TPD+\L
D)(R+D^TPD)^{-1}(B^TP+D^TPC+D^T\L)X,X\ran\big]dt\Big\} \\ [2mm]
\dsp={1\over2}\dbE\Big\{\int_0^\t\(\lan Q_1X_1,X_1\ran+\lan
R_1u,u\ran\)dt+\lan G_1(\t)X_1(\t),X_1(\t)\ran\\ [3mm]
\dsp\q+\int_\t^T\(\lan(R+D^TPD)^{-1}(B^TP+D^TPC+D^T\L)X,(B^TP+D^TPC+D^T\L)X\ran \\ [3mm]
\dsp\q+2\lan(B^TP+D^TPC+D\L)X,u\ran+\lan(R+D^TPD)u,u\ran\)dt\Big\} \\ [3mm]
\dsp=\frac{1}{2}\dbE\Big\{\int_0^\t\(\lan Q_1X_1,X_1\ran+\lan R_1u,u\ran\)dt+\lan G_1(\t)X_1(\t),X_1(\t)\ran \\ [3mm]
\dsp\q+\int_\t^T\left|(R+D^TPD)^{1\over2}\[u+(R+D^TPD)^{-1}(B^TP+D^TPC+D^T\L)X\]\right|^2dt\Big\}.\ea$$
Therefore, we have
$$\ba{ll}
\dsp J^{\t}(x_1;u(\cd))= \frac{1}{2}\dbE\int_0^\t\(\lan
Q_1X_1,X_1\ran+\lan
R_1u,u\ran\)dt+\lan[G_1(\t)+K(\t)^TP(\t)K(\t)]X_1(\t),X_1(\t)\ran \\
[3mm]
\dsp\qq\qq\qq\q+\int_\t^T\left|(R+D^TPD)^{1\over2}\[u+(R+D^TPD)^{-1}(B^TP+D^TPC+D^T\L)X\]\right|^2dt\Big\}.\ea$$
Next, applying It\^o formula to $\lan P_1(\cd)X_1(\cd),X_1(\cd)\ran$
on the interval $[0,\t]$, similar to the above, we have
$$\ba{ll}
\dsp J^{\t}(x_1;u(\cd))={1\over2}\dbE\Big\{\lan P_1(0)x_1,x_1\ran \\ [2mm]
\dsp\q+\int_0^\t\left|(R_1+D_1^TP_1D_1)^{1\over2}\[u+(R_1+D_1^TP_1D_1)^{-1}(B_1^TP_1+D_1^TP_1C_1+D_1^T\L_1)
X_1\]\right|^2dt\\ [3mm]
\dsp\q+\int_\t^T\left|(R+D^TPD)^{1\over2}\[u+(R+D^TPD)^{-1}(B^TP+D^TPC+D^T\L)X\]\right|^2dt\Big\}\\
\dsp\equiv{1\over2}\dbE\Big\{\lan P_1(0)x_1,x_1\ran+\int_0^\t|(R_1+D_1^TP_1D_1)^{-1}(u+\Psi_1X_1)|^2dt+\int_\t^T|(R+D^TPD)^{-1}(u+\Psi X)|^2dt\Big\}.\ea$$
Then our conclusion follows.
\endpf

\section{The Equivalence of Control Problems}

Under {\bf(H3.1)}, Problem (RLQ) admits a unique optimal pair which
can be represented by (\ref{bar u}). In this section, we will
establish some connection between Problem (MLQ) and Problem (RLQ).
To this end, we denote the solution to Riccati equation
(\ref{P_1-Riccati0}) on $[0,\t]$ by $(P_1^\t(\cd),\L_1^\t(\cd))$,
emphasizing its dependence on $\t$ via the terminal condition. It is
clear that $\t\mapsto P^\t_1(s)$ is continuous. Therefore, the
following problem makes sense:\\

\fbox{\textbf{Problem (OT)}} For given $x_1\in\dbR^{n_1}$, find
a $\bar\t\in\cT[0,T]$ such that
\bel{bar tau}\lan P_1^{\bar\t}(0)x_1,x_1\ran=\inf_{\t\in\cT[0,T]}\lan P_1^\t(0)x_1,x_1\ran.\ee
Any $\bar\t \in\cT[0,T]$ satisfying the above is called an optimal
time of Problem (OT). Note that in the case $n_1=1$, if $\bar\t$ is
an optimal time for some $x_1\in\dbR\setminus\{0\}$, then it
is an optimal time for all $x_1\in\dbR\setminus\{0\}$. On the
other hand, since there might not be some kind of monotonicity of
the map $\t\mapsto P_1^{\t}(s)$, optimal time $\bar\t$ may not be unique.
Now we establish the equivalence between Problem (MLQ) and Problems
(RLQ) and (OT).

\ms

\bf Theorem 4.1. \sl Suppose \textbf{(H3.1)} holds and
$(\bar{u}(\cdot), \bar{\tau})$ is an optimal control pair of Problem
(MLQ). Then the optimal control $\bar{u}(\cdot)$ can be represented by
\bel{4.2}\bar u(t)=\begin{cases}-\Psi(t)\bar X(t),& t \in [\bar\t,T],\\
-\Psi_1^{\bar\t}(t)\bar X_1(t), & t \in [0,\bar\t),
\end{cases}\ee
where
\bel{4.3}\left\{\3n\ba{ll}
\ns\ds\Psi(t)=\[R(t)+D(t)^TP(t)D(t)\]^{-1}\[B(t)^TP(t)+D(t)^TP(t)C(t)+D(t)^T\L(t)\],\\
\ns\ds\Psi^{\bar\t}_1(t)=\[R_1(t)+D_1(t)^TP^{\bar\t}_1(t)D_1(t)\]^{-1}\[B_1(t)^TP^{\bar\t}_1(t)+D(t)^TP^{\bar\t}_1(t)C_1(t)+D_1(t)^T\L_1^{\bar\t}(t)\].\ea\right.\ee
Moreover, the optimal value of the cost functional for Problem (MLQ) is given by
\bel{J(bar u)}J(x_1;\bar u(\cd),\bar\t)={1\over2}\lan P^{\bar\t}_1(0)x_1,x_1\ran,\ee
with $\bar\t$ being an optimal time for Problem (OT) to Riccati
equation $(\ref{P-Riccati0})$ and $(\ref{P_1-Riccati0})$.

\ms

\it Proof. \rm Suppose $(\bar{u}(\cdot), \bar{\tau})$ is an optimal
pair of Problem (MLQ), that is
$$J(x_1;\bar u(\cd),\bar\t)=\inf_{(u(\cd),\t)\in\cU[0,T]\times\cT[0,T]}J(x_1;u(\cd),\t).$$
We fix $\bar\t$. Then the above implies
$$J^{\bar\t}(x_1;\bar u(\cd))=\inf_{u(\cd)\in\cU[0,T]}J^{\bar\t}(x_1;u(\cd)).$$
Then it follows from Theorem $3.1$ that $\bar u(\cd)$ admits
representation (\ref{4.2}) with $\Psi(\cd)$ and
$\Psi_1^{\bar\t}(\cd)$ given by (\ref{4.3}), and
$$J^{\bar\t}(x_1;\bar u(\cd))={1\over2}\lan P^{\bar\t}_1(0)x_1,x_1\ran.$$
Next, for any stopping time $\t \in \mathcal{T}[0,T]$, we can
construct a control $\widetilde{u}\in \mathcal{U}[0, T]$ satisfying
$$\wt u(t)=\begin{cases}-\Psi(t)\widetilde{X}(t),& t \in [\tau, T],\\
-\Psi_1^\t(t)\widetilde{X}_1(t), & t \in [0,\tau),\end{cases}$$
where $\widetilde{X}(\cdot), \widetilde{X}_1(\cdot)$ are
defined in similar way to $\bar{X}(\cdot), \bar{X}_1(\cdot)$ but replacing
$\bar{\t}$ by $\t$. Following the similar arguments, we can prove
$$J^\t(x_1;\wt u(\cd))=\inf_{u(\cd)\in\cU[0,T]}J^\t(x_1;u(\cd))=\frac{1}{2}\lan {P}_1^{\tau}(0)x_1,x_1\ran.$$
Then
$$\ba{ll}
\ns\ds{1\over2}\lan P_1^{\bar\t}(0)x_1, x_1\ran=J^{\bar\t}(x_1;\bar u(\cd))=J(x_1;\bar u(\cd),\bar\t)\\
\ns\ds=\inf_{(u(\cd),\t)\in\cU[0,T]\times\cT[0,T]}J(x_1; u(\cd),\t)\le J(x_1;\wt u(\cd),\t)={1\over2}\lan P_1^\t(0)x_1,x_1\ran.\ea$$
That is, $\bar{\tau}$ solves the Problem (OT). Hence
the results.
\endpf

\ms

Next, we consider the case of deterministic coefficients. More
precisely, we introduce the following assumption.

\ms

\textbf{(H3.1$'$)} The following holds:
\bel{}\left\{\ba{ll}
\ns\ds A_1(\cd),C_1(\cd)\in L^\infty(0,T;\dbR^{n_1\times n_1}),\qq A(\cd),C(\cd)\in L^\infty(0,T;\dbR^{n\times n}),\\
\ns\ds B_1(\cd),D_1(\cd)\in L^\infty(0,T;\dbR^{n_1\times m}),\qq B(\cd),D(\cd)\in L^\infty(0,T;\dbR^{n\times m}),\\
\ns\ds Q_1(\cd)\in L^\infty(0,T;\dbS^{n_1}),\qq Q(\cd)\in L^\infty(0,T;\dbS^n),\\
\ns\ds R(\cd)\in L^\infty(0,T;\dbS^m),\qq R_1(\cd)\in L^\infty(0,T;\dbS^m),\\
\ns\ds G_1(\cd)\in C([0,T];\dbS^{n_1}),\q G\in\dbS^n,\q K(\cd)\in
C([0,T];\dbR^{n\times n_1}).\ea\right.\ee
Moreover, for some $\d>0$,
\bel{}\left\{\ba{ll}
\ns\ds R_1(t)\ge\d I_{n_1},\q R(t)\ge\d I_n,\\
\ns\ds Q_1(t),G_1(t)\ge0,\q Q(t)\ge0,\q G\ge0,\ea\right.\qq
t\in[0,T].\ee

We have the following result.

\ms

\bf Proposition 4.2. \sl Let \textbf{(H3.1$'$)} hold. Let $P(\cd)$
solves
\bel{PP-Riccati}\left\{\ba{ll}
\ns\ds\dot{P}+PA+A^TP+C^TPC+Q-(PB+C^TPD)(R+D^TPD)^{-1}\big(B^TP+D^TPC\big)=0,\\
\ns\ds\qq\qq\qq\qq\qq\qq\qq\qq\qq\qq t\in[0,T],\\
\ns\ds P(T)=G,\qq R+D^TPD>0,\ea\right.\ee
and $P_1^r(\cd)$ solves
\bel{P_11-Riccati}\left\{\ba{ll}
\ns\ds
\dot P^r_1+P^r_1A_1+A_1^TP^r_1+C_1^TP^r_1C_1+Q_1\\
\ns\ds\qq\qq-(P^r_1B_1+C_1^TP^r_1D_1)(R_1+D_1^TP^r_1D_1)^{-1}
\big(B_1^TP^r_1+D_1^TP^r_1C_1\big)=0,\q t\in[0,r], \\
\ns\ds P^r_1(r)=K(r)^{T}P(r)K(r)+G_1(r),\qq R_1+D_1^T\bar
P^r_1D_1>0,\ea\right.\ee
with $r\in(0,T)$. Then Problem (MLQ) has an optimal control pair
$(\bar u(\cd),\bar r)\in\cU[0,T]\times[0,T]$ where
\bel{}\bar u(t)=\left\{\2n\ba{ll}
-\Psi(t)\bar X(t),\qq\q t\in[\bar r,T],\\
-\Psi^{\bar{r}}_1(t)\bar X_1(t),\qq t\in[0,\bar r),\ea\right.\ee
with $\bar r$ being deterministic and
\bel{}\left\{\ba{ll}
\ns\ds\Psi(t)=\[R(t)+D(t)^TP(t)D(t)\]^{-1}\[B(t)^TP(t)+D(t)^TP(t)C(t)\],\\
\ns\ds\Psi_1^{\bar r}(t)=\[R_1(t)+D_1(t)^TP_1^{\bar
r}(t)D_1(t)\]^{-1}\[B_1(t)^TP_1^{\bar r}(t)+D(t)^TP_1^{\bar
r}(t)C_1(t)\].\ea\right.\ee
Moreover,
\bel{bar r}\lan P^{\bar r}_1(0)x_1,x_1\ran=\inf_{r\in[0,T]}\lan
P^r_1(0)x_1,x_1\ran.\ee

\ms

\it Proof. \rm In the current case, Riccati equation
(\ref{P-Riccati0}) becomes (\ref{PP-Riccati}). On the other hand, if
we still denote $(P^\t_1(\cd),\L^\t_1(\cd))$ to be the adapted
solution to the Riccati equation (\ref{P_1-Riccati0}), then when
$\t=r\in(0,T)$ is deterministic, one has $\L^r_1(\cd)=0$ and
$P^r_1(\cd)$ satisfies deterministic Riccati equation
(\ref{P_11-Riccati}). It is clear that if $\bar r\in[0,T]$ satisfies
(\ref{bar r}), then
$${1\over2}\lan P_1^r(0)x_1,x_1\ran=\inf_{\t\in\cT[0,T]}{1\over2}
\lan P_1^\t(0)x_1,x_1\ran=V(x_1).$$
This proves our result. \endpf

\ms

For the convenience below, we state the following problem.

\ms

\fbox{\textbf{Problem (DOT)}} Find $\bar r\in[0,T]$ such that
\bel{}\lan P^{\bar r}_1(0)x_1,x_1\ran=\inf_{r\in[0,T]}\lan
P^r_1(0)x_1,x_1\ran.\ee

\ms

By Proposition $4.2$, we see that the optimal time $\bar r$ in Problem
(DOT) solves Problem (MLQ). Moreover, in general, the optimal time
for Problem (DOT) depends on the Riccati equation $P^{\bar
r}_1(\cdot)$ and the initial condition $x_1$. The following example
makes this clear.

\ms

\bf Example 4.3. \rm Let $n_1=2$, $n_2=1$, and let
$$\left\{\3n\ba{ll}
\ns\ds A\1n=\1n\left(\ba{ccc}0&0&0\\0&0&0\\0&0&a\ea\right),~
B\1n=\1n\left(\ba{c}0\\0\\
1\ea\right),~C\1n=\1n\left(\ba{ccc}0&0&0\\0&0&0\\0&0&0\ea\right),~
D\1n=\1n\left(\ba{c}0\\0\\
0\ea\right),\\ [5mm]
\ns\ds A_1\1n=\1n\left(\ba{cc}0&1\\0&0\ea\right),~
B_1\1n=\1n\left(\ba{c}0\\
1\ea\right),~C_1\1n=\1n\left(\ba{cc}0&0\\0&0\ea\right),~
D_1\1n=\1n\left(\ba{c}0\\
0\ea\right),\\
[3mm]
\ns\ds Q\1n=\1n\left(\ba{ccc}0&0&0\\0&0&0\\0&0&0\ea\right),~
R\1n=\1n1,~G=\left(\ba{ccc}0&0&0\\0&0&0\\0&0&g\ea\right),\\
\ns\ds
Q_1\1n=\1n\left(\ba{cc}0&0\\0&0\ea\right),~R_1\1n=\1n1,~G_1=\left(\ba{cc}g_1&0\\0&0\ea\right),~
K\1n=\1n\left(\ba{cc}1&0\\0&1\\1&0\ea\right),\ea\right.$$
with $a\in\dbR$, $g,g_1\in(0,\infty)$. Then for any $r\in(0,T)$, we
have the state equation
$$\left\{\ba{ll}
\ns\ds\dot X_1^1(t)=X_1^2(t),\qq\qq t\in[0,r),\\
\ns\ds\dot X_1^2(t)=u(t),\qq\qq t\in[0,r),\\
\ns\ds\dot X_2(t)=aX_2(t)+u(t),\qq t\in[r,T],\\
\ns\ds X_1(0)=x_1,\qq X_2(r)=X_1^1(r).\ea\right.$$
The cost functional reads
$$\ba{ll}
\ns\ds
J(x_1;u(\cd),r)={1\over2}\Big\{\int_0^r|u(t)|^2dt+g_1|X_1^1(r)|^2
+\int_r^T|u(t)|^2dt+g|X_2(T)|^2\Big\}. \ea$$
In this case, we have
$$P(t)=\left(\ba{ccc}0&0&0\\0&0&0\\0&0& P_2(t)\ea\right),\qq t\in[r,T],$$
with $P_2(\cd)$ solves the following Riccati equation:
$$\left\{\ba{ll}
\ns\ds\dot P_2(t)+2aP_2(t)-P_2(t)^2=0,\qq t\in[r,T],\\
\ns\ds P_2(T)=g.\ea\right.$$
This equation admits a unique solution $P_2(\cd)$. We claim that,
\bel{P_2-rep}P_2(t)=\left\{\ba{ll}
\ns\ds{2age^{2a(T-t)}\over g(e^{2a(T-t)}-1)+2a},\qq t\in[r,T],\qq
a\ne0,\\
\ns\ds{g\over1+g(T-t)},\qq\qq\q t\in[r,T],\qq a=0.\ea\right.\ee
Note that if $a>0$, then
$$g(e^{2a(T-t)}-1)+2a\ge2a>0,$$
and if $a<0$, then
$$g(e^{2a(T-t)}-1)+2a\le2a<0.$$
Hence,
$$P_2(t)>0,\qq t\in[r,T].$$

\ms

Next, the Riccati equation for $P_1(\cd)$ reads
$$\left\{\ba{ll}
\ns\ds\dot P_1(t)+P_1(t)A_1+A_1^TP_1(t)-P_1(t)MP_1(t)=0,\qq
t\in[0,r],\\
\ns\ds P_1(r)=\bar G_1(r),\ea\right.$$
where
$$M=B_1B_1^T=\left(\ba{cc}0&0\\0&1\ea\right),\q
\bar G_1(r)=K^TP(r)K+G_1=\left(\ba{cc}g_1+P_2(r)&0\\0&0\ea\right).$$
We now solve the above Riccati equation by a method found in
\cite{MY}. To this end, let
$$\bar P_1(t)=P_1(t)-\bar G_1(r),\qq t\in[0,r].$$
Then (suppressing $t$ and $r$)
$$\ba{ll}
\ns\ds\dot{\bar P}_1=\dot P_1=-\[P_1A_1+A_1^TP_1-P_1MP_1\]\\
\ns\ds\q~=-\[(\bar P_1+\bar G_1)A_1+A_1^T(\bar P_1+\bar G_1)-(\bar
P_1+\bar G_1)M(\bar P_1+\bar G_1)\]\\
\ns\ds\q~=-\[\bar P_1(A_1-M\bar G_1)+(A_1-M\bar G_1)^T\bar P_1-\bar
P_1M\bar P_1+\bar G_1A_1+A_1^T\bar G_1-\bar G_1M\bar G_1\].\ea$$
Note that
$$M\bar G_1=\left(\ba{cc}0&0\\0&1\ea\right)\left(\ba{cc}g_1+P_2(r)&0\\0&0\ea\right)=0,$$
$$\bar
G_1A_1=\left(\ba{cc}g_1+P_2(r)&0\\0&0\ea\right)\left(\ba{cc}0&1\\0&0\ea\right)
=\left(\ba{cc}0&g_1+P_2(r)\\0&0\ea\right).$$
Hence, $\bar P_1(\cd)$ should be the solution to the following
Riccati equation:
$$\left\{\ba{ll}
\ns\ds\dot{\bar P}_1(t)+\bar P_1(t)A_1+A_1^T\bar P_1(t)-\bar
P_1(t)M\bar P_1(t)+\bar g_1J=0,\qq t\in[0,r],\\
\ns\ds\bar P_1(r)=0,\ea\right.$$
where
$$J=\left(\ba{cc}0&1\\1&0\ea\right),\qq\bar g_1=g_1+P_2(r).$$
Next, we let
$$\cA=\left(\ba{cc}A_1&-M\\-\bar g_1J&-A_1^T\ea\right)\equiv\left(\ba{cccc}0&1&0&0\\
                                                                          0&0&0&-1\\
                                                                          0&-\bar g_1&0&0\\
                                                                          -\bar g_1&0&-1&0\ea\right).$$
Then according to \cite{MY}, we have the following representation of
the solution $P_1(\cd)$:
$$\bar P_1(t)=-\[(0,I)e^{\cA(r-t)}\left(\ba{c}0\\
I\ea\right)\]^{-1}(0,I)e^{\cA(r-t)}\left(\ba{c}I\\0\ea\right),\qq
t\in[0,r],$$
as long as the involved inverse exists. We now calculate
$e^{\cA(r-t)}$. Direct computation show that
$$\ba{ll}
\ns\ds\cA^2=\left(\ba{cccc}0&0&0&-1\\
                       \bar g_1&0&1&0\\
                       0&0&0&\bar g_1\\
                       0&0&0&0\ea\right),\q
%
\cA^3=\left(\ba{cccc}\bar g_1&0&1&0\\
                       0&0&0&0\\
                       -\bar g_1^2&0&-\bar g_1&0\\
                       0&0&0&0\ea\right),\q
\cA^4=0.\ea$$
Hence,
$$\ba{ll}
\ns\ds e^{\cA(r-t)}=I+(r-t)\cA+{(r-t)^2\over2}\cA^2+{(r-t)^3\over6}\cA^3\\
\ns\ds=\left(\ba{cccc}1+{\bar
g_1(r-t)^3\over6}&(r-t)&{(r-t)^3\over6}&-{(r-t)^2\over2}\\
{\bar g_1(r-t)^2\over2}&1&{(r-t)^2\over2}&-(r-t)\\
-{\bar g_1^2(r-t)^3\over6}&-\bar g_1(r-t)&1-{\bar
g_1(r-t)^3\over6}&{\bar g_1(r-t)^2\over2}\\
-\bar
g_1(r-t)&0&-(r-t)&1\ea\right)\equiv\left(\ba{cc}\F_{11}&\F_{12}\\
\F_{21}&\F_{22}\ea\right).\ea$$
Then
$$\ba{ll}
\ns\ds(0,I)e^{\cA(r-t)}\left(\ba{c}0\\
I\ea\right)=(0,I)\left(\ba{cc}
\F_{11}&\F_{12}\\ \F_{21}&\F_{22}\ea\right)\left(\ba{c}0\\
I\ea\right)=\F_{22}=\left(\ba{cc}1-{\bar g_1\over6}(r-t)^3&{\bar
g_1\over2}(r-t)^2\\-(r-t)&1\ea\right).\ea$$
Since
$$\det\left(\ba{cc}1-{\bar g_1\over6}(r-t)^3&{\bar
g_1\over2}(r-t)^2\\-(r-t)&1\ea\right)=1+{\bar
g_1(r-t)^3\over3}>0,\qq\forall t\in[0,r],$$
we have
$$\left(\ba{cc}1-{\bar g_1\over6}(r-t)^3&{\bar
g_1\over2}(r-t)^2\\-(r-t)&1\ea\right)^{-1}={3\over3+\bar
g_1(r-t)^3}\left(\ba{cc}1&-{\bar g_1\over2}(r-t)^2\\(r-t)&1-{\bar
g_1\over6}(r-t)^3\ea\right).$$
On the other hand,
$$(0,I)e^{\cA(r-t)}\left(\ba{c}I\\0\ea\right)=(0,I)\left(\ba{cc}\F_{11}&\F_{12}\\ \F_{21}&\F_{22}\ea\right)
\left(\ba{c}I\\0\ea\right)=\F_{21}=\left(\ba{cc}-{\bar
g_1^2(r-t)^3\over6}&-\bar g_1(r-t)\\-\bar g_1(r-t)&0\ea\right).$$
Hence,
$$\ba{ll}
\ns\ds\bar P_1(t)=-\[(0,I)e^{\cA(r-t)}\left(\ba{c}0\\
I\ea\right)\]^{-1}(0,I)e^{\cA(r-t)}\left(\ba{c}I\\0\ea\right)=-\F_{22}^{-1}\F_{21}\\
\ns\ds=-{3\over3+\bar g_1(r-t)^3}\left(\ba{cc}1&-{\bar
g_1\over2}(r-t)^2\\(r-t)&1-{\bar
g_1\over6}(r-t)^3\ea\right)\left(\ba{cc}-{\bar
g_1^2(r-t)^3\over6}&-\bar g_1(r-t)\\-\bar
g_1(r-t)&0\ea\right)\\
\ns\ds=-{3\over3+\bar g_1(r-t)^3}\left(\ba{cc}{\bar
g_1^2\over3}(r-t)^3&-\bar g_1(r-t)\\-\bar g_1(r-t)&-\bar
g_1(r-t)^2\ea\right)={\bar g_1(r-t)\over3+\bar
g_1(r-t)^3}\left(\ba{cc}-\bar g_1(r-t)^2&3\\3&3(r-t)\ea\right).\ea$$
Consequently,
$$P_1(t)=\bar P_1(t)+\bar{G}_1(r)={3[g_1+P_2(r)]\over3+[g_1+P_2(r)](r-t)^3}\left(\ba{cc}1&r-t\\
r-t&(r-t)^2\ea\right),\qq t\in[0,r],$$
which is positive definite on $[0,r)$. Hence,
$$\ba{ll}
\ns\ds V(r,x_1)=\lan
P^r_1(0)x_1,x_1\ran={3[g_1+P_2(r)]\over3+[g_1+P_2(r)]r^3}\[(x_1^1)^2+2rx_1^1x_1^2+r^2(x_1^2)^2\]\\
\ns\ds\qq\qq={3[g_1+P_2(r)]\over3+[g_1+P_2(r)]r^3}\big[x_1^1+rx_1^2\big]^2.\ea$$
Clearly, for different $x_1\in\dbR^2$, the optimal $\bar r$ will be
different in general.

\bs

\section{One-Dimensional Cases with Constant Coefficients}

In this section, we make the following assumption
\bel{5.1}\left\{\ba{ll}
\ns\ds n_1=n_2=m=1,\\ 
\ns\ds A(t)=\left(\ba{cc}0&0\\0&A_2\ea\right),\q B=\left(\ba{c}0\\
B_2\ea\right),\q C(t)=\left(\ba{cc}0&0\\
0&C_2\ea\right),\q D(\cd)=\left(\ba{c}0\\ D_2\ea\right),\\
\ns\ds Q(t)=\left(\ba{cc}0&0\\0&Q_2\ea\right),\q
R(t)=R_2,\q G=\left(\ba{cc}0&0\\
0&G_2\ea\right),\q K(\cdot)=\left(\ba{c}1\\
K\ea\right),\\ [2mm]
\ns\ds A_1(t)=A_1,\q B_1(t)=B_1,\q C_1(t)=C_1,\q
D_1(t)=D_1,\\
\ns\ds Q_1(t)=Q_1,\q R_1(t)=R_1,\q G_1(t)=G_1,\ea\right.\ee
where $A_i,B_i,C_i,D_i,Q_i,R_i,G_i,K$ are all constants ($i=1,2$), and
\bel{5.2}R_1,R_2>0,\q Q_1,Q_2,G_1,G_2\ge0,\q K\ne0.\ee
Then the controlled system becomes
\bel{state}\left\{\2n\ba{ll}
\ns\ds dX_1(t)=\big[A_1X_1(t)+B_1u(t)\big]dt+\big[C_1X_1(t)+D_1u(t)\big]dW(t),\qq t\in[0,\t),\\
\ns\ds dX_2(t)=\big[A_2X_2(t)+B_2u(t)\big]dt+\big[C_2X_2(t)+D_2u(t)\big]dW(t),\qq t\in[\t,T], \\
\ns\ds X_1(0)=x_1,\qq X_2(\tau)=KX_1(\t-0), \ea\right. \ee
and the cost functional is
\bel{cost2}\ba{ll}
\ns\ds J(x_1;u(\cd),\t)={1\over2}\dbE\Big\{\int_0^\t\big[Q_1X_1(t)^2
+R_1u(t)^2\big]dt+G_1X_1(\t)^2\\
\ns\ds\qq\qq\qq\qq\q+\int_\t^T\big[Q_2X_2(t)^2+R_2u(t)^2\big]dt+
G_2X_2(T)^2\Big\}.\ea\ee
Thus, the first component $X_1(\cd)$ of the state process will be
completely terminated from $\t$ on. In this case, we
have
$$P(t)=\left(\ba{cc}0&0\\0&P_2(t)\ea\right),\qq t\in[r,T],$$
where $P_2(\cd)$ satisfies
\bel{P-Riccati}\left\{\2n\ba{ll}
\ns\ds\dot
P_2(t)+(2A_2+C_2^2)P_2(t)+Q_2-{(B_2+C_2D_2)^2P_2(t)^2\over
R_2+D_2^2P_2(t)}=0,\qq t\in[0,T],\\
\ns\ds P_2(T)=G_2,\qq R_2+D_2^2P_2(t)>0,\ea\right.\ee
and $P^r_1(\cd)$ satisfies
\bel{P_1-Riccati}\left\{\2n\ba{ll}
\ns\ds\dot
P^r_1(t)+(2A_1+C_1^2)P^r_1(t)+Q_1-{(B_1+C_1D_1)^2P_1^r(t)^2\over
R_1+D_1^2P^r_1(t)}=0,\qq t\in[0,r],\\
\ns\ds P_1^r(r)=K^2P_2(r)+G_1,\qq R_1+D_1^2P_1(r)>0.\ea\right.\ee
If we denote
$$F_i(P)=(2A_i+C_i^2)P+Q_i-{(B_i+C_iD_i)^2P^2\over R_i+D_i^2P},\qq
i=1,2.$$
Then (\ref{P-Riccati}) and (\ref{P_1-Riccati}) can be written as
\bel{P-Riccait2}\left\{\ba{ll}
\ns\ds\dot P_2(t)+F_2\big(P_2(t)\big)=0,\qq t\in[0,T],\\
\ns\ds P_2(T)=G_2,\ea\right.\ee
and
\bel{P_1-Riccati2}\left\{\ba{ll}
\ns\ds \dot P_1^r(t)+F_1\big(P^r_1(t)\big)=0,\qq t\in[0,r],\\
\ns\ds P_1^r(r)=K^2P_2(r)+G_1.\ea\right.\ee
Note that since $n_1=1$, the optimal time $\bar r$ of Problem (DOT)
is independent of the initial state $x_1$. Thus, the optimal time
$\bar r$ satisfies
\bel{5.7}P_1^{\bar r}(0)=\inf_{r\in[0,T]}P_1^r(0).\ee
The following gives a necessary condition for $\bar r$.

\ms

\bf Proposition 5.2. \sl Let $(\ref{5.1})$--$(\ref{5.2})$ hold. Then
the optimal time $\bar{r}$ to Problem {\rm(DOT)} satisfies the
following condition:
\bel{}F_1\big(K^2P_2(\bar r)+G_1\big) -K^2F_2\big(P_2(\bar
r)\big)\left\{\ba{ll}
\ns\ds\ge0,\qq\qq\bar r=0,\\
\ns\ds=0,\qq\qq\bar r\in(0,T),\\
\ns\ds\le0,\qq\qq\bar r=T.\ea\right.\ee

\ms

\it Proof. \rm We first claim that
$$\Pi^r(t)={d\over dr}P^r_1(t),\qq t\in[0,T],$$
exists for any $r\in(0,T)$ and $\Pi^r(\cd)$ solves
\bel{Pi}\left\{\ba{ll}
\ns\ds\dot\Pi^r(t)+F_1'\big(P^r_1(t)\big)\Pi^r(t)=0,\qq t\in[0,r],\\
\ns\ds\Pi^r(r)=F_1\big(K^2P_2(r)+G_1\big)-K^2F_2\big(P_2(r)\big).\ea\right.\ee
To see this, let $r\in(0,T)$ and $\e>0$ small so that
$r\pm\e\in(0,T)$. Consider the following: for any $t\in[0,r]$,
$$\ba{ll}
\ns\ds
P_1^{r+\e}(t)-P_1^r(t)=K^2\[P_2(r+\e)-P_2(r)\]+\int_r^{r+\e}F_1\big(P^{r+\e}_1(s)\big)ds\\
\ns\ds\qq\qq\qq\qq\qq+\int_t^r\[F_1\big(P_1^{r+\e}(s)
\big)-F_1\big(P_1^r(s)\big)\]ds.\ea$$
Then by a standard argument, we have the existence of the following
limit:
$$\Pi^r_+(t)\equiv\lim_{\e\da0}{P^{r+\e}_1(t)
-P^r_1(t)\over\e},\qq t\in[0,r],$$
and
$$\ba{ll}
\ns\ds\Pi^r_+(t)=K^2\dot P_2(r)+F_1\big(P_1^r(r)\big)+\int_t^rF'\big(P_1^r(s)
\big)\Pi^r_+(s)ds\\
\ns\ds\qq=-K^2F_2\big(P_2(r)\big)+F_1\big(K^2P_2(r)+G_1\big)+\int_t^rF'\big
(P_1^r(s)\big)\Pi^r_+(s)ds,\qq t\in[0,r].\ea$$
On the other hand,
$$\ba{ll}
\ns\ds
P_1^{r-\e}(t)-P_1^r(t)=K^2\[P_2(r-\e)-P_2(r)\]-\int_{r-\e}^rF_1\big(P_1^r(s)\big)ds\\
\ns\ds\qq\qq\qq\qq\qq+\int_t^{r-\e}\[F_1\big(P_1^{r-\e}(s)\big)-F_1\big(P_1^r(s)\big)\]ds.\ea$$
Again, by a standard argument, we have the existence of the
following limit:
$$\Pi^r_-(t)\equiv\lim_{\e\da0}{P^{r-\e}_1(t)-P^r_1(t)\over-\e},\qq t\in[0,r],$$
and
$$\ba{ll}
\ns\ds\Pi^r_-(t)=K^2\dot P_2(r)+F_1\big(P_1^r(r)\big)+\int_t^rF_1'\big(P_1^{r}(s)\big)\Pi^r_-(s)ds\\
\ns\ds\qq\q=-K^2F_2\big(P_2(r)\big)+F_1\big(K^2P_2(\bar r)+G_1\big)+\int_t^rF_1'\big(P^r_1(s)\big)\Pi^r_-(s)ds.\ea$$
Thus, $r\mapsto P_1^r(t)$ is differentiable with derivative
$${d\over dr}P^r_1(t)=\Pi^r(t)\equiv\Pi^r_\pm(t),\qq t\in[0,T],$$
and $\Pi(\cd)$ satisfies (\ref{Pi}). Clearly,
$$\Pi^r(t)=e^{\int_t^rF_1'\big(P_1^r(s)\big)ds}\[F_1\big(K^2P_2(r)+G_1\big)
-K^2F_2\big(P_2(r)\big)\],\qq t\in[0,r].$$
Now, since $r\mapsto P^r_1(0)$ attains a minimum at $\bar
r\in(0,T)$, we have
$$\Pi^{\bar r}(0)=0,$$
which leads to our conclusion for $\bar r\in(0,T)$. Now, if $\bar
r=0$, then
$$\Pi^0(0)=\lim_{\e\da0}{P_1^\e(0)-P_1^0(0)\over\e}\ge0.$$
Finally, if $\bar r=T$, then
$$\Pi^T(0)=\lim_{\e\da0}{P_1^{T-\e}(0)-P_1^T(0)\over-\e}\le0.$$
This completes the proof. \endpf

\ms

If the optimal time $\bar r$ is either 0 or $T$, our
Problem (MLQ) becomes less interesting. Therefore, the optimal time
$\bar r$ is said to be {\it non-trivial} if $\bar r\in(0,T)$. From
the above, one has the following corollary.

\ms

\bf Corollary 5.2. \sl Let $(\ref{5.1})$--$(\ref{5.2})$ hold. Then
\bel{0<r,r<T}\left\{\ba{ll}
\ns\ds F_1(K^2G_2+G_1)-K^2F_2(G_2)>0\qq\Ra\qq\bar r<T,\\ [3mm]
\ns\ds F_1(K^2P_2(0)+G_1)-K^2F_2(P_2(0))<0\qq\Ra\qq\bar
r>0.\ea\right.\ee
Therefore, $\bar r$ is non-trivial if
that
\bel{0<r<T}F_1(K^2G_2+G_1)-K^2F_2(G_2)>0>F_1(K^2P_2(0)+G_1)-K^2F_2(P_2(0)).\ee

\ms

\rm

Note that in principle, conditions in (\ref{0<r<T}) are checkable.
Let us now look at some special cases for which we can say something
about the optimal time $\bar r$. Let
\bel{}D_2=G_1=0,\q R_2=K=1,\q B_2\ne0.\ee
In this case, we observe the following:
$$\ba{ll}
\ns\ds F_1(K^2P+G_1)-K^2F_2(P)=F_1(P)-F_2(P)\\
\ns\ds\q=(2A_1+C_1^2)P+Q_1-{(B_1+C_1D_1)^2P^2\over
R_1+D_1^2P}-\[(2A_2+C_2^2)P+Q_2-B^2_2P^2\]\\
\ns\ds\q=\[2(A_1-A_2)+C_1^2-C_2^2\]P+Q_1-Q_2+B_2^2P^2
-{(B_1+C_1D_1)^2P^2\over R_1+D_1^2P}\equiv{\Th(P)\over
R_1+D_1^2P},\ea$$
where
$$\ba{ll}
\ns\ds\Th(P)=D_1^2B_2^2P^3+\[\(2(A_1-A_2)+C_1^2-C_2^2\)D_1^2+R_1B_2^2
-(B_1+C_1D_1)^2\]P^2\\
\ns\ds\qq\qq+\[\(2(A_1-A_2)+C_1^2-C_2^2\)R_1+(Q_1-Q_2)D_1^2\]P+(Q_1-Q_2)R_1.\ea$$
Then (\ref{0<r<T}) is implied by
\bel{G_2>P_2}\Th(G_2)>0>\Th\big(P_2(0)\big).\ee
On the other hand, in the current case, the Riccati equation for
$P_2(\cd)$ becomes
\bel{}\left\{\2n\ba{ll}
\ns\ds\dot P_2(t)+(2A_2+C_2^2)P_2(t)+Q_2-B_2^2P_2(t)^2=0,\qq t\in[0,T],\\
\ns\ds P_2(T)=G_2.\ea\right.\ee
Denote
$$\l_\pm={2A_2+C_2^2\pm\sqrt{(2A_2+C_2^2)^2+4B_2^2Q_2}\over2B_2^2}.$$
Then we may rewrite the above Riccati equation as
\bel{}\left\{\2n\ba{ll}
\ns\ds\dot P_2(t)-B_2^2[P_2(t)-\l_+][P_2(t)-\l_-]=0,\qq t\in[0,T],\\
\ns\ds P_2(T)=G_2.\ea\right.\ee
If
$$G_2=\l_\pm,$$
then the unique solution $P_2(\cd)$ is given by
$$P_2(t)\equiv G_2,\qq t\in[0,T],$$
for which, (\ref{G_2>P_2}) cannot be true. Therefore, in what
follows, we assume that
\bel{G_2 ne l}G_2\ne\l_\pm.\ee
Then the solution $P_2(\cd)$ is given by
$$P_2(t)={\l_+(G_2-\l_-)e^{B_2^2(\l_+-\l_-)(T-t)}-\l_-(G_2-\l_+)\over(G_2-\l_-)e^{B_2^2(\l_+-\l_-)(T-t)}
-(G_2-\l_+)},\qq t\in[0,T].$$
We claim that
$$P_2(0)={\l_+(G_2-\l_-)e^{B_2^2(\l_+-\l_-)T}-\l_-(G_2-\l_+)
\over(G_2-\l_-)e^{B_2^2(\l_+-\l_-)T}-(G_2-\l_+)}\ne G_2.$$
In fact, it is easy to see that there is no $t_0\in(0,T)$ such that
$$B_2^2[P_2(t_0)-\l_+][P_2(t_0)-\l_-]=\dot P_2(t_0)=0.$$
Otherwise, by the uniqueness of solutions to ODEs, we must have
$$P_2(t)\equiv\left\{\ba{ll}\l_+,\qq t\in[0,T],\qq\hb{if }P_2(t_0)-\l_+=0,\\
\ns\ds\l_-,\qq t\in[0,T],\qq\hb{if }P_2(t_0)-\l_-=0,\ea\right.$$
both of which contradict (\ref{G_2 ne l}). Actually, when (\ref{G_2
ne l}) holds, by observing the sign of $\dot P_2(\cd)$, one has the
following:
\bel{}P_2(0)\left\{\ba{ll}
\ns\ds>G_2,\qq\hb{if }G_2\in(\l_-,\l_+),\\
\ns\ds<G_2,\qq\hb{if }G_2\notin[\l_-,\l_+].\ea\right.\ee
In particular, if $Q_2=0$, then
$$\l_+={2A_2+C_2^2\over B_2^2},\qq\l_-=0.$$
Consequently,
$$P_2(0)={G_2(2A_2+C_2^2)e^{(2A_2+C_2^2)T}
\over G_2B_2^2e^{(2A_2+C_2^2)T}+2A_2+C_2^2-G_2B_2^2}\left\{\ba{ll}
\ns\ds>G_2,\qq2A_2+C_2^2>G_2B_2^2,\\ [3mm]
\ns\ds<G_2,\qq2A_2+C_2^2<G_2B_2^2.\ea\right.$$
Further, let $D_1=0$ and $R_1=1$. Then
$$\ba{ll}
\ns\ds\Th(P)=(B_2^2-B_1^2)P^2+\(2A_1+C_1^2-(2A_2+C_2^2)\)P+Q_1.\ea$$
Let
$$B_2^2<B_1^2.$$
Then $\Th(\cd)$ has a unique positive root:
$$P_+={(2A_1+C_1^2)-(2A_2+C_2^2)+\sqrt{\big[(2A_1+C_1^2)-(2A_2+C_2^2)\big]^2
+4(B_1^2-B_2^2)Q_1}\over2(B_1^2-B_2^2)}.$$
Hence, if
$$2A_2+C_2^2>G_2B_2^2,$$
and
$$0<G_2<P_+<P_2(0),$$
then (\ref{G_2>P_2}) holds and $\bar r$ is non-trivial.

\ms

It is clear that many other cases for which $\bar r$ is non-trivial
can be discussed in the similar fashion. However, we prefer to omit
the details here.

\ms

\end{document}